\documentstyle[12pt]{article}

\begin{document}
\newcommand{\ol }{\overline}
\newcommand{\ul }{\underline }
\newcommand{\ra }{\rightarrow }
\newcommand{\lra }{\longrightarrow }
\newcommand{\ga }{\gamma }
\newcommand{\st }{\stackrel }
\newcommand{\scr }{\scriptsize }
\title{{\Large\bf Polynilpotent Multipliers of Finitely Generated Abelian
Groups \footnote{This research was in part supported by a grant
from IPM (No.82200029).}}}
\author{{\bf Behrooz Mashayekhy and Mohsen Parvizi\footnote{E-mail address: mashaf@math.um.ac.ir \ \ \ \
\
parvizi@math.um.ac.ir}} \\
Department of Mathematics, Ferdowsi University of Mashhad,\\
P.O.Box 1159-91775, Mashhad, Iran \\ and\\
Institute for Studies in Theoretical Physics and Mathematics, \\
P.O.Box 5746-19395, Tehran, Iran}
\date{ }
\maketitle \vspace{-0.4cm}
 \hspace{-0.7cm}\line(1,0){390}\\
\vspace{-0.4cm}
\hspace{-0.2cm}{\bf Abstract}\\

{\footnotesize In this paper, we present an explicit formula for
the Baer invariant of a finitely generated abelian group with
respect to the variety of polynilpotent groups of class row
$(c_1,\ldots ,c_t)$, ${\cal N}_{c_1,\ldots ,c_t}$. In particular,
one can obtain an explicit structure of the $\ell$-solvable
multiplier ( the Baer invariant with respect to the vaiety of
solvable groups of length at most $\ell\geq 1,\ {\cal S}_{\ell}$.)
of a finitely generated abelian group.\\

\hspace{-0.65cm}{\it Mathematics Subject Classification(2000)}: 20E34, 20E10, 20F18.}\\
\hspace{-0.7cm}\line(1,0){390}\\

\hspace{-0.7cm}{\bf 1. Introduction and preliminaries}\\

I.Schur [13], in 1907, found a formula for the Schur multiplier of
a direct product of two finite groups as follows:
$$ M(A\times B)\cong M(A)\oplus M(B)\oplus A_{ab}\otimes B_{ab}.$$
One of the important corollaries of the above fact is an explicit
formula for the Schur multiplier of a finite abelian group $G\cong
{\bf Z}_{n_1}\oplus {\bf Z}_{n_2}\oplus \ldots ,{\bf Z}_{n_k}$,
where $n_{i+1}\mid n_i$ for all $1\leq i\leq k-1$, as follows:
$$M(G)\cong {\bf Z}_{n_2}\oplus {\bf Z}_{n_3}^{(2)}\oplus \ldots \oplus {\bf Z}_{n_k}^{(k-1)},$$
where ${\bf Z}_n^{(m)}$ denotes the direct sum of $m$ copies of
the cyclic group ${\bf Z}_n$ (see [10]).

In 1997, the first author, in a joint paper [11], succeeded to
generalize the above formula for the Baer invariant of a finite
abelian group $G \cong {\bf Z}_{n_1} \oplus {\bf Z}_{n_2} \oplus
\ldots ,{\bf Z}_{n_k}$, where $n_{i+1}\mid n_i$ for all $1\leq
i\leq k-1$, with respect to the variety of nilpotent groups of
class at most $c\geq 1$, ${\cal N}_c$, as follows:

$${\cal N}_cM(G)\cong {\bf Z}_{n_2}^{(b_2)}\oplus {\bf Z}_{n_3}^{(b_3-b_2)}\oplus \ldots \oplus
{\bf Z}_{n_k}^{(b_k-b_{k-1})},$$ where $b_i$ is the number of
basic commutators of weight $c+1$ on $i$ letters (see [4]).

${\cal N}_cM(G)$ is also called the $c$-nilpotent multiplier of
$G$ (see [3]). Note that, by a similar method of the paper [11],
we can obtain the structure of the $c$-nilpotent multiplier of a
finitely generated abelian group as the following theorem.\\

\hspace{-0.8cm} {\bf Theorem 1.1.} Let $G\cong {\bf Z}^{(m)}\oplus
{\bf Z}_{n_1}\oplus {\bf Z}_{n_2}\oplus \ldots {\bf Z}_{n_k}$ be a
finitely generated abelian group, where $n_{i+1}\mid n_i$ for all
$1\leq i\leq k-1$, then
$${\cal N}_cM(G)\cong {\bf Z}^{(b_m)}\oplus {\bf Z}_{n_1}^{(b_{m+1}-b_m)}\oplus {\bf Z}_{n_2}^{(b_{m+2}-b_{m+1})}\oplus \ldots
\oplus {\bf Z}_{n_k}^{(b_{m+k}-b_{m+k-1})}\ ,$$ where $b_i$ is the
number of basic commutators of weight $c+1$ on $i$ letters and
$b_0=b_1=0$.

 Now, in this paper, we intend to generalize the above theorem to
obtain an explicit formula for ${\cal N}_{c_1,\ldots ,c_t}M(G)$,
the Baer invariant of $G$ with respect to the variety of
polynilpotent groups of class row $(c_1,\ldots ,c_t)$, ${\cal
N}_{c_1,\ldots ,c_t}$, where $G$ is a finitely generated abelian
group. We also call ${\cal N}_{c_1,\ldots ,c_t}M(G)$, a
polynilpotent multiplier of $G$. As an immediate consequence, one
can obtain an explicit formula for the $\ell$-solvable multiplier
of
$G$, ${\cal S}_{\ell}M(G)$.\\

\hspace{-0.6cm}{\bf Definition 1.2.} Let $G$ be any group with a
free presentation $G\cong F/R$, where $F$ is a free group. Then,
after R. Baer [1], the Baer invariant of $G$ with respect to a
variety of groups ${\cal V}$, denoted by ${\cal V}M(G)$, is
defined to be
$${\cal V}M(G)=\frac{R\cap V(F)}{[RV^*F]}\ ,$$
where $V$ is the set of words of the variety ${\cal V}$, $V(F)$ is
the verbal subgroup of $F$ with respect to ${\cal V}$ and
$$[RV^*F]=<v(f_1,\ldots ,f_{i-1},f_ir,f_{i+1},\ldots,f_n)v(f_1,\ldots,f_i,
\ldots,f_n)^{-1}\mid $$
$$r\in R, 1\leq i\leq n, v\in V ,f_i\in F, n\in {\bf N}>.$$

In special case, if ${\cal V}$ is the variety of abelian groups,
${\cal A}$ , then the Baer invariant of $G$ will be
$$\frac{R\cap F'}{[R,F]},$$
which, following Hopf [7], is isomorphic to the second cohomology
group of $G$, $H_2(G,C^*)$, in finite case and also is isomorphic
to the well-known notion the Schur multiplier of $G$, denoted by
$M(G)$. The multiplier $M(G)$ arose in Schur's work [12] of 1904
on projective representation of a group, and has subsequently
found a variety of other applications. The survey article of
Wiegold [14] and the books of Beyl and Tappe [2] and Karpilovsky
[10] form a fairly comprehensive account of $M(G)$.

If ${\cal V}$ is the variety of nilpotent groups of class at most
$c\geq1$, ${\cal N}_c$, then the Baer invariant of $G$ with
respect to ${\cal N}_c$ will be
$${\cal N}_cM(G)=\frac{R\cap \gamma_{c+1}(F)}{[R,\ _cF]},$$
where $\gamma_{c+1}(F)$ is the $(c+1)$-st term of the lower
central series of $F$ and $[R,\ _1F]=[R,F], [R,\ _cF]=[[R,\
_{c-1}F],F]$, inductively.

If ${\cal V}$ is the variety of solvable groups of length at most
$\ell\geq1$, ${\cal S}_{\ell}$, then the Baer invariant of $G$
with respect to ${\cal S}_{\ell}$ will be
$${\cal S}_{\ell}M(G)=\frac{R\cap
\delta^l(F)}{[R,F,\delta^1(F),\ldots,\delta^{l-1}(F)]},$$ where
$\delta^i(F)$ is the $i$-th derived subgroup of $F$.See
[8,corollary 2.10] for the equality
$[RS_{l}^*F]=[R,F,\delta^1(F),\ldots,\delta^{l-1}(F)]$.

 In a very more general case, let $\cal V$ be the variety of
 polynilpotent groups of class row $(c_1,\ldots ,c_t)$, ${\cal
 N}_{c_1,\ldots ,c_t}$, then the Baer invariant of a group $G$
 with respect to this variety is as follows:
 $${\cal N}_{c_1,\ldots ,c_t}M(G)\cong \frac{R\cap
 \ga_{c_t+1}\circ\ldots\circ\ga_{c_1+1}(F)}{[R,\ _{c_1}F,\
 _{c_2}\ga_{c_1+1}(F),\ldots ,
 \ _{c_t}\ga_{c_{t-1}+1}\circ\ldots\circ \ga_{c_1+1}(F)]}\ ,$$
 where
 $\ga_{c_t+1}\circ\ldots\circ\ga_{c_1+1}(F)=\ga_{c_t+1}(\ga_{c_{t-1}+1}
 (\ldots(\ga_{c_1+1}(F))\ldots
 ))$ are the terms of iterated lower central series of $F$. See [
 6, Corollary 6.14] for the equality
 $$[R{\cal N}_{c_1,\ldots ,c_t}^*F]=[R,\ _{c_1}F,
 \ _{c_2}\ga_{c_1+1}(F),\ldots ,
\ _{c_t}\ga_{c_{t-1}+1}\circ \ldots \circ \ga_{c_1+1}(F)].$$

In the following, we are going to mention some definitions and
notations of T.C. Hurley and M.A. Ward [9], which are vital in our
investigation.\\
\newpage
\hspace{-0.8cm}
 {\bf Definition and Notation 1.3.} Commutators are written $[a,b]=a^{-1}b^{-1}ab$ and the usual
 convention for left-normed commutators is used,
 $[a,b,c]=[[a,b],c],\ [a,b,c,d]=[[[a,b],c],d]$ and so on,
 including the trivial case $[a]=a$.

 {\it Basic commutators} are defined in the usual way. If $X$ is a fully
 ordered independent subset of a free group, the basic commutators
 on $X$ are defined inductively over their weight as follows:\\
$(i)$ All the members of $X$ are basic commutators on $X$ of
weight one on $X$.\\
$(ii)$ Assuming that $n>1$ and that the basic commutators of
weight less than $n$ on $X$ have been defined and ordered.\\
$(iii)$ A commutator $[a,b]$ is a basic commutator of weight $n$
on $X$ if  $wt(a)+wt(b)=n,\ a<b$, and if $b=[b_1,b_2]$, then
$b_2\leq a$. The ordering of basic commutators is then extended to
include those of weight $n$ in any way such that those of weight
less than $n$ precede those of weight $n$. The natural way to
define the order on basic commutators of the same weight is
lexicographically, $[b_1,a_1]<[b_2,a_2]$ if $b_1<b_2$ or if
$b_1=b_2$ and $a_1<a_2$.

 A word of the form
$$[c,a_1,a_2,\ldots ,a_p,b_1^{\beta_1},b_2^{\beta_2},\ldots,b_q^{\beta_q}]$$
is a {\it ``standard invertator"} will be meant to imply that the
$\beta_i$'s are $\pm 1$, $c>a_1\leq a_2\leq \ldots\leq a_p\leq
b_1\leq b_2\leq \ldots\leq b_q$ and if $b_i=b_j$ then
$\beta_i=\beta_j$ for all $i,j$. Whenever this terminology is used
it will be accomplished by a statement of what set $X$, the $a_i$
and the $b_j$ are chosen from and this will be always be a set
which is known to be fully ordered in some way. Restrictions on
the values of $p$ and $q$ will be given, the value $p=0$ and $q=0$
being permissible so that we may, when we wish, specify standard
invertators of the forms $[c,a_1,\ldots,a_p]$ or
$[c,b_1^{\beta_1},b_2^{\beta_2},\ldots,b_q^{\beta_q}].$

 Let $F$ be a free group on alphabet $X$ and $m$ and $n$ be integers.
Then\\
$(i)$ $A_{m,n}$ denotes the set of all basic commutators on $X$ of
weight exactly $n$ and of the form $[c,a_1,\ldots,a_p]$, where $b$
and the $a_i$ are all basic commutators on $X$ of weight less than
$m$.\\
$(ii)$ $B_{m,n}$ denotes the set of all standard invertators on
$X$ of the form
$$[b,a_1,a_2,\ldots,a_p,a_{p+1}^{\alpha_{p+1}},\ldots,a_q^{\alpha_q}],$$
where $0\leq p<q$, $b$ and the $a_i$ are basic commutators on $X$
of weight less than $m$, $$wt([b,a_1,a_2,\ldots ,a_p])<n\leq
wt([b,a_1,a_2,\ldots ,a_p,a_{p+1}^{\alpha_{p+1}}])$$ and
$b=[b_1,b_2]$ implies $b_2\leq a_1$. Note that $[b,a_1,a_2,\ldots
,a_p]\in A_{m,r}$, where $r$ is the weight of this commutator and
$r<n$. Also, observe that $A_{m,m}$ is just the set of all basic
commutators of weight $m$ on $X$.\\

\hspace{-0.6cm}{\bf Theorem 1.4} (P.Hall [4,5]). Let
$F=<x_1,x_2,\ldots ,x_d>$ be a free group, then
$$ \frac {\ga_n(F)}{\ga_{n+i}(F)} \ \ , \ \ \ \  1\leq i\leq n$$
is the free abelian group freely generated by the basic
commutators of weights
$n,n+1,\ldots ,n+i-1$ on the letters $\{x_1,\ldots ,x_d\}.$\\

\hspace{-0.6cm}{\bf Theorem 1.5} (Witt Formula [4]). The number of
basic commutators of weight $n$ on $d$ generators is  given by the
following formula:
$$ \chi_n(d)=\frac {1}{n} \sum_{m|n}^{} \mu (m)d^{n/m}$$
where $\mu (m)$ is the {\it Mobious function}, and defined to be
   \[ \mu (m)=\left \{ \begin{array}{ll}
      1 & {\rm if}\ m=1, \\ 0 & {\rm if}\ m=p_1^{\alpha_1}\ldots
p_k^{\alpha_k}\ \ \exists \alpha_i>1, \\ (-1)^s & {\rm if}\
m=p_1\ldots p_s,
\end{array} \right.  \]

 The following important theorem presents interesting free
generating sets for the terms of the lower central series of a
free group which will be used several times in this paper.\\

\hspace{-0.6cm}{\bf Theorem 1.6} (T.C.Hurley and M.A.Ward 1981).
Let $F$ be a free group, freely generated by some fully ordered
set $X$, and let $m$ and $n$ be integers satisfying $2\leq m\leq
n$. Then the members of the set
$$A_{m,m}\cup A_{m,m+1}\cup \ldots \cup A_{m,n-1}\cup B_{m,n}$$
are distinct as written, so that in particular this is a disjoint
union, and the set freely generates $\gamma_m(F)$.\\

\hspace{-0.6cm}{\bf Proof.} See [9, Theorem 2.2].\\

\hspace{-0.6cm}{\bf Corollary 1.7.} Let $F$ be a free group freely
generated by some fully ordered set $X$. Then
$\ga_{c_2+1}(\ga_{c_1+1}(F))$ is freely generated by
$$\hat{A}_{c_2+1,c_2+1}\cup \hat{B}_{c_2+1,c_2+2},$$
where $\hat{A}_{c_2+1,c_2+1}$ is the set of all basic commutators
of weight $c_2+1$ on the set $$Y=A_{c_1+1,c_1+1}\cup
B_{c_1+1,c_1+2},$$
 and $ \hat{B}_{c_2+1,c_2+2}$ is the set of all
standard invertators on $Y$ of the form
$$[b,a_1,a_2,\ldots,a_p,a_{p+1}^{\alpha_{p+1}},\ldots,a_q^{\alpha_q}],$$
where $0\leq p<q$, $b$ and the $a_i$ are basic commutators on $Y$
of weight less than $c_2+1$, $$wt([b,a_1,a_2,\ldots
,a_p])<c_2+2\leq wt([b,a_1,a_2,\ldots
,a_p,a_{p+1}^{\alpha_{p+1}}])$$ and $b=[b_1,b_2]$ implies $b_2\leq
a_1$.\\

\hspace{-0.6cm}{\bf Proof.} Using Theorem 1.6, $\ga_{c_1+1}(F)$ is
freely generated by $A_{c_1+1,c_1+1}\cup B_{c_1+1,c_1+2}$, when
putting $m=c_1+1, n=c_1+2$. Now we can suppose
$\overline{F}=\ga_{c_1+1}(F)$ is a free group, freely generated by
fully ordered set $Y=A_{c_1+1,c_1+1}\cup B_{c_1+1,c_1+2}$.
Applying Theorem 1.6 again for $\ga_{c_2+1}(\overline{F})$ and
$m=c_2+1, n=c_2+2$, the result holds. $\Box$\\

 As an immediate consequence we have the following corollary.\\

\hspace{-0.6cm}{\bf Corollary 1.8.} Let $F$ be a free group freely
generated by some fully ordered set $X$. Then the second derived
subgroup of $F$, $\delta^2(F)=F''$, is freely generated by
$$\hat{A}_{2,2}\cup \hat{B}_{2,3},$$
where $\hat{A}_{2,2}$ is the set of all basic commutators of
weight 2 on the set $A_{2,2}\cup B_{2,3}$, and $\hat{B}_{2,3}$ is
the set of all standard invertators on $A_{2,2}\cup B_{2,3}$ of
the form $[b,a_1,a_2^{\alpha _2},\ldots,a_q^{\alpha_q}]$, where
$b,a_i\in A_{2,2}\cup B_{2,3}$.\\

\hspace{-0.6cm}{\bf 2. The Main Results}

 In this section, first, we concentrate on the calculation of the
Baer invariant of a finitely generated abelian group with respect
to the variety of metabelian groups, $i.e.$ solvable groups of
length 2, ${\cal S}_2$.

 Let ${\bf Z}_{r_i}=<x_i \mid x_i^{r_i}> , 1\leq i\leq t $, be cyclic
groups of order $r_i\geq 0$, and let
$$0\longrightarrow R_i=<x_i^{r_i}>\longrightarrow F_i=<x_i>\longrightarrow {\bf Z}_{r_i}
\longrightarrow 0\ ,$$ be a free presentation of ${\bf Z}_{r_i}$.
Also, suppose $G\cong \oplus{\sum_{i=1}^{t}} {\bf Z}_{r_i}$ is the
direct sum of the cyclic groups ${\bf Z}_{r_i}$. Then
$$0\longrightarrow R\longrightarrow F \longrightarrow G \longrightarrow 0$$
is a free presentation of $G$, where $F={\prod_{i=1}^{*^t}}
F_i=<x_1,\ldots,x_t>$ is the free product of $F_i$'s, and
$R=\prod_{i=1}^t R_i\gamma_2(F)$. Therefore, the metabelian
multiplier of $G$ is as follows:
$${\cal S}_2M(G)\cong \frac{R\cap \delta^2(F)}{[R,F,\delta^1(F)]}=
\frac{F''}{[R,F,F']}\ \ \ \  (since \ \ \ F'\leq R).$$

 Now, the following theorem presents an explicit structure for the
metabelian multiplier of a finitely generated abelian group.\\

\hspace{-0.6cm}{\bf Theorem 2.1.} With the above notation and
assumption, let $G\cong {\bf Z}^{(m)}\oplus {\bf Z}_{n_1}\oplus
\ldots \oplus {\bf Z}_{n_k}$ be a finitely generated abelian
group, where $n_{i+1}\mid n_i$ for all $1\leq i\leq k-1$. Then the
following isomorphism holds:
$${\cal S}_2M(G)\cong {\bf Z}^{(d_m)}\oplus {\bf Z}_{n_1}^{(d_{m+1}-d_m)}\oplus
\ldots \oplus {\bf Z}_{n_k}^{(d_{m+k}-d_{m+k-1})},$$ where
$d_i=\chi_2(\chi_2(i))$, and $\chi_2(i)$ is the number of
basic commutators of weight 2 on $i$ letters.\\

\hspace{-0.6cm}{\bf Proof.} With the previous notation, put
$t=m+k$, $r_1=r_2=\ldots=r_m=0$, $r_{m+j}=n_j$ , $1\leq j\leq k$.
Then ${\bf Z}_{r_1}\cong \ldots \cong {\bf Z}_{r_m}\cong {\bf Z}$,
${\bf Z}_{r_{m+j}}\cong {\bf Z}_{n_j}$, $G\cong \oplus \sum
_{i=1}^{m+k} {\bf Z}_{r_i}$, and
$${\cal S}_2M(G)\cong \frac{F''}{[R,F,F']},$$
where $F$ is the free group on the set
$X=\{x_1,\ldots,x_m,x_{m+1},\ldots,x_{m+k}\}.$

By corollary 1.8 $F''$ is a free group with the basis
$\hat{A}_{2,2}\cup \hat{B}_{2,3}$. Put $L$ the normal closure of
those elements of the basis $F''$, $\hat{A}_{2,2}\cup
\hat{B}_{2,3}$, of weight, as commutators on the set
$X=\{x_1,\ldots,x_m,x_{m+1},\ldots,x_{m+k}\},$ greater than 4 in
$F''$. In other words
$$L=<w\in \hat{A}_{2,2}\cup \hat{B}_{2,3} \mid w\notin
\{u\in \hat{A}_{2,2}\cup \hat{B}_{2,3}\mid u\ \ is\ \ of\ \ the\ \
form\ \ [[x_{i_1},x_{i_2}],[x_{i_3},x_{i_4}]]\}>^{F''}.$$ It is
easy to see that $F''/L$ is a free group freely generated by the
following set
$$Y=\{wL \mid w\in \hat{A}_{2,2}\cup \hat{B}_{2,3}
 \ \ and\ w\ \ is\ \ of\ \ the\ \ form\ \ [[x_{i_1},x_{i_2}],[x_{i_3},x_{i_4}]]\}.$$
Therefore
$$\frac{F''/L}{(F''/L)'}\cong
\frac{F''}{LF'''}$$ is a free abelian group with the basis
$\overline{Y}=\{wLF'''\mid wL\in Y\}.$ Since ${\cal S}_2M(G)\cong
F''/[R,F,F']$ is abelian, so $F'''\leq [R,F,F']$. Thus, we have
$${\cal S}_2M(G)\cong \frac{F''/LF'''}{[R,F,F']/LF'''}.$$

Now we are going to describe explicitly the bases of the free
abelian group $F''/LF'''$ and its subgroup $[R,F,F']/LF'''$ in
order to find the structure of the metabelian multiplier of $G$,
${\cal S}_2M(G)$. According to the basis $\overline{Y}$ of the
free abelian group $F''/LF'''$, it is easy to see that
$$\overline{Y}=C_0\cup C_1\cup \ldots \cup C_k\ ,$$
where
$$C_0=\{wLF'''\in\overline{Y}\mid
w=[[x_{i_1},x_{i_2}],[x_{i_3},x_{i_4}]], 1\leq i_1,i_2,i_3,i_4\leq
m\},$$ and for all $1\leq \lambda \leq k$
$$C_{\lambda}=\{wLF'''\in\overline{Y}\mid
w=[[x_{i_1},x_{i_2}],[x_{i_3},x_{i_4}]], 1\leq i_1,i_2,i_3,i_4\leq
m+\lambda $$
$$,\exists \ 1\leq j\leq 4,\ \ s.t.\ i_j=m+\lambda\}.$$

In order to find an appropriate basis for the free abelian group
$[R,F,F']/LF'''$, first we claim that $\gamma_5(F)\cap F''\leq
LF'''$ (*), since, let $u\in F''$, using the basis $\overline{Y}$
of the free abelian group $F''/LF'''$, we have
$$uLF'''=w_{i_1}^{\epsilon_1}\ldots w_{i_t}^{\epsilon_t} LF''',$$
where $w_{i_1}LF''',\ldots ,w_{i_t}LF''' \in \overline{Y}$, and
$\epsilon_1,\ldots,\epsilon_t\in {\bf Z}$. Clearly $LF'''\leq
\gamma_5(F)$, so, if $u\in \gamma_5(F)$, then we have
$w_{i_1}^{\epsilon_1}\ldots w_{i_t}^{\epsilon_t}\in \gamma_5(F)$.
It is easy to see that $w_{i_1}^{\epsilon_1},\ldots
,w_{i_t}^{\epsilon_t}$ are basic commutators of weight 4 on $X$.
By Theorem 1.4 $\gamma_4(F)/\gamma_5(F)$ is the free abelian group
with basis of all basic commutators of weight 4 on $X$. Thus we
have $\epsilon_1=\ldots=\epsilon_t=0$, and hence $u\in LF'''$. As
an immediate consequence we have $[F',F,F']\leq LF'''$. Note that
$R=(\prod_{i=1}^{m+k}R_i)F'$, where $R_i=<x_i^0>=1$, for all
$1\leq i\leq m$, and $R_{m+j}=<x_{m+j}^{n_j}>$, for all $1\leq
j\leq k$, so
$$\frac{[R,F,F']}{LF'''}=\frac{\prod_{j=1}^{k}[R_{m+j},F,F']LF'''}{LF'''}\ .$$
Using the above equality and the congruence
$$[[x_{i_1}^{\alpha_1},x_{i_2}^{\alpha_2}],[x_{i_3}^{\alpha_3},x_{i_4}^{\alpha_4}]]
\equiv[[x_{i_1},x_{i_2}],[x_{i_3},x_{i_4}]]^{{\alpha_1}{\alpha_2}{\alpha_3}{\alpha_4}}
  \ \ (mod \ \ LF''')\ ,$$
for all $\alpha_1,\alpha_2,\alpha_3,\alpha_4\in {\bf Z}$, (By
(*)), it is routine to check that the free abelian group
$[R,F,F']/LF'''$ has the following basis
$$D_1\cup D_2\cup\ldots\cup D_k\ ,$$
where $D_{\lambda}=\{w^{n_{\lambda}}LF''' \mid wLF'''\in
C_{\lambda}$, $1\leq \lambda \leq k\}.$

Using the form of the elements $C_{\lambda}$ and the number of
basic commutators of weight 2 on $i$ letters, $\chi_2(i)$, one can
easily see that $\mid C_0\mid =\chi_2(\chi_2(m))$, and $\mid
C_{\lambda}\mid=\chi_2(\chi_2(m+\lambda))-\chi_2(\chi_2(m+\lambda-1)).$
Hence the result holds.$\Box$\\

Now, trying to generalize the proof of the above theorem, which is
the basic idea of the paper, we are going to present an explicit
formula for the polynilpotent multiplier of a finitely generated
abelian group with respect to the variety ${\cal N}_{c_1,\ldots
,c_t}$. Because of applying an iterative method and avoiding
complicacy for the reader, first, we state and prove the beginning
step of the method for the variety ${\cal N}_{c_1,c_2}$ in the
following
theorem.\\

\hspace{-0.6cm}{\bf Theorem 2.2.} Let ${\cal N}_{c_1,c_2}$ be the
polynilpotent variety of class row $(c_1,c_2)$ and $G\cong {\bf
Z}^{(m)}\oplus {\bf Z}_{n_1}\oplus\ldots\oplus {\bf Z}_{n_k}$ be a
finitely generated abelian group, where $n_{i+1}\mid n_i$ for all
$1\leq\i\leq k-1$. Then the following isomorphism holds:
$${\cal N}_{c_1,c_2}M(G)\cong {\bf Z}^{(e_m)}\oplus {\bf Z}_{n_1}^{(e_{m+1}-e_m)}\oplus\ldots\oplus
{\bf Z}_{n_k}^{(e_{m+k}-e_{m+k-1})},$$ where
$e_i=\chi_{c_2+1}(\chi_{c_1+1}(i))$ for all $m\leq i\leq
m+k$.\\

\hspace{-0.6cm}{\bf Proof.} By the notation of the Theorem 2.1 we
have
$${\cal N}_{c_1,c_2}M(G)=\frac{\gamma_{c_2+1}(\gamma_{c_2+1}(F))}
{[R,\ _{c_1}F,\ _{c_2}\gamma_{c_1+1}(F)]},$$ where $F$ is the free
group on the set $X=\{x_1,\ldots,x_m,x_{m+1},\ldots,x_{m+k}\}$. By
considering the basis of the free group
$\gamma_{c_2+1}(\gamma_{c_2+1}(F))$ presented in corollary 1.7, we
put
$$L=<w\in \hat{A}_{c_2+1,c_2+1}\cup \hat{B}_{c_2+1,c_2+2}\mid \
w\notin E>^{\gamma_{c_2+1}(\gamma_{c_1+1}(F))},$$ where $E$ is the
set of all basic commutators of weight exactly $c_2+1$ on the set
of all basic commutators of weight exactly $c_1+1$ on the set $X$.

Clearly $\gamma_{c_2+1}(\gamma_{c_1+1}(F))/L$ is free on the set
$$Y=\{wL\mid \ w\in \hat{A}_{c_2+1,c_2+1}\cup
\hat{B}_{c_2+1,c_2+2}\ and \ w\in E\}$$ and
$\gamma_{c_2+1}(\gamma_{c_1+1}(F))/L\gamma_{c_2+1}(\gamma_{c_1+1}(F))$
is free abelian with the basis
$\bar{Y}=\{wL\gamma_{c_2+1}(\gamma_{c_1+1}(F))\mid \ wL\in Y\}$.
Considering the form of the elements of $L$ and noticing to the
abelian group ${\cal N}_{c_1,c_2}M(G)$, we have
$$L\gamma_{c_2+1}(\gamma_{c_1+1}(F))\leq [R,\ _{c_1}F,\ _{c_2}\gamma_{c_1+1}(F)].$$

Thus the following isomorphism holds:
$${\cal N}_{c_1,c_2}M(G)\cong \frac{\gamma_{c_2+1}(\gamma_{c_1+1}(F))/
L\gamma_2(\gamma_{c_2+1}(\gamma_{c_1+1}(F)))}{[R,\ _{c_1}F,\
_{c_2}\gamma_{c_1+1}(F)]/
L\gamma_2(\gamma_{c_2+1}(\gamma_{c_1+1}(F)))}.$$ By Theorem 1.4
$\gamma_{c_1+c_2+c_1c_2+1}(F)/\gamma_{c_1+c_2+c_1c_2+2}(F)$ is the
free abelian group with the basis of all basic commutators of
weight $c_1+c_2+c_1c_2+1$ on $X$. Using the above fact and the
basis $\bar{Y}$ of the free abelian group\\
$\gamma_{c_2+1}(\gamma_{c_1+1}(F))/L\gamma_{c_2+1}(\gamma_{c_1+1}(F))$
we can conclude the following inclusion:
$$\gamma_{c_1+c_2+c_1c_2+2}(F)\cap \gamma_{c_2+1}(\gamma_{c_1+1}(F))
\leq L\gamma_2(\gamma_{c_2+1}(\gamma_{c_1+1}(F))).$$

Now, it is easy to see that $\bar{Y}=C_0\cup C_1\cup\ldots\cup
C_k$ is a basis for the free abelian group
$\gamma_{c_2+1}(\gamma_{c_1+1}(F))/L\gamma_2(\gamma_{c_2+1}(\gamma_{c_1+1}(F)))$
and $ D_1\cup D_2\cup\ldots\cup D_k$ is a basis for the free
abelian group
$$\frac{[R,\ _{c_1}F,\ _{c_2}\gamma_{c_1+1}(F)]}{L\gamma_2(\gamma_{c_2+1}(\gamma_{c_1+1}(F)))},$$
where
$$C_0=\{wL\gamma_2(\gamma_{c_2+1}(\gamma_{c_1+1}(F)))\in \bar{Y}\mid \ w\in E$$
\begin{center}
 and $w$ is a commutator on letters $x_1,\ldots,x_m\},$
\end{center}
and for $1\leq \lambda \leq k$:\\
\hspace{-0.7cm}$C_{\lambda}=\{wL\gamma_2(\gamma_{c_2+1}(\gamma_{c_1+1}(F)))\in
\bar{Y}\mid \ w\in E $  and $w$ is a commutator on letters
 $x_1,\ldots,x_m,x_{m+1},\ldots,
 x_{m+\lambda}$
 such that the letter $x_{m+\lambda}$ does appear in $w\},$\\
\hspace{-0.7cm}$D_{\lambda}=\{w^{n_{\lambda}}L\gamma_2(\gamma_{c_2+1}(\gamma_{c_1+1}(F)))\in
\bar{Y}\mid \ wL\gamma_2(\gamma_{c_2+1}(\gamma_{c_1+1}(F)))\in
C_{\lambda}\}.$\\
\hspace{-0.7cm} Note that using the form of the elements of
$C_{\lambda}$ and the number of basic commutators, we can conclude
that $\mid C_0\mid=\chi_{c_2+1}(\chi_{c_1+1}(m))$ and $\mid
C_{\lambda}\mid=\chi_{c_2+1}(\chi_{c_1+1}(m+\lambda))-\chi_{c_2+1}(\chi_{c_1+1}(m+\lambda-1)).$
Hence the result holds. $\Box$\\

Now, we are ready to state and prove the main result of the paper
in general case.\\

\hspace{-0.6cm}{\bf Theorem 2.3.} Let ${\cal
N}_{c_1,c_2,\ldots,c_t}$ be the polynilpotent variety of class row
$(c_1,c_2,\ldots,c_t)$ and $G\cong {\bf Z}^{(m)}\oplus {\bf
Z}_{n_1}\oplus\ldots\oplus {\bf Z}_{n_k}$ be a finitely generated
abelian group, where $n_{i+1}\mid n_i$ for all $1\leq\i\leq k-1$.
Then an explicit structure of the polynilpotent multiplier of $G$
is as follows.
$${\cal N}_{c_1,c_2,\ldots,c_t}M(G)\cong {\bf Z}^{(f_m)}\oplus {\bf Z}_{n_1}^{(f_{m+1}-f_m)}\oplus\ldots\oplus
{\bf Z}_{n_k}^{(f_{m+k}-f_{m+k-1})},$$ where
$f_i=\chi_{c_t+1}(\chi_{c_{t-1}+1}(\ldots(\chi_{c_1+1}(i))\ldots))$
for all $m\leq i \leq m+k$.\\

\hspace{-0.6cm}{\bf Proof.} Let $F$ be the free group on the set
$X=\{x_1,\ldots,x_m,x_{m+1},\ldots,x_{m+k}\}$. Then by previous
notation, we have
$${\cal N}_{c_1,c_2,\ldots,c_t}M(G)=\frac{\gamma_{c_t+1}(\gamma_{c_{t-1}+1}(\ldots(\gamma_{c_1+1}(F))\ldots))}{[R,\ _{c_1}F,\ _{c_2}\gamma_{c_1+1}(F),\ldots,\ _{c_t}\gamma_{c_t+1}(\gamma_{c_{t-1}+1}(\ldots(\gamma_{c_1+1}(F))\ldots))]}.$$
We define $\rho_t(F),E_t,X_t$ inductively on $t$ as follows:\\
$\rho_1(F)=\gamma_{c_1+1}(F),\
\rho_i(F)=\gamma_{c_i+1}(\rho_{i-1}(F);\\
\ E_1=X,\ E_i=$the set of all basic commutators of weight $c_i+1$ on the set $E_{i-1}$;\\
$X_1=A_{c_1+1,c_1+1}\cup B_{c_1+1,c_1+2},\ X_i={\hat
A}_{c_i+1,c_i+1}\cup {\hat B}_{c_i+1,c_i+2},$ where ${\hat
A}_{c_i+1,c_i+1}$ is the set of all basic commutators of weight
$c_i+1$ on the set $X_{i-1}$, and ${\hat B}_{c_i+1,c_i+2}$ is the
set of all standard invertators on $X_{i_1}$ of the form
$$[b,a_1,\ldots,a_p,a_{p+1}^{\alpha_{p+1}},\ldots,a_q^{\alpha_q}],$$
where $0\leq p<q$, $b$ and the $a_i$ are basic commutators on
$X_{i-1}$ of weight less than $c_i+1$ ,
$wt([b,a_1,\ldots,a_p])<c_i+2\leq
wt([b,a_1,\ldots,a_p,a_{p+1}^{\alpha_{p+1}}])$ and $b=[b_1,b_2]$
implies $b_2\leq a_1$.

Using Theorem 1.6 and induction on $t$, it is easy to see that\\
$\gamma_{c_t+1}(\gamma_{c_{t-1}+1}(\ldots(\gamma_{c_1+1}(F))\ldots))=\rho_t(F)$
is freely generated by $X_t$. Now, putting
$$L_t=<w\in X_t \mid \ w\notin E_t>^{\rho_t(F)},$$
one can easily see that $\rho_t(F)/L_t$ is free on the set
$Y_t=\{wL_t \mid \ w\in E_t\}$ and
$\rho_t(F)/L_t\gamma_2(\rho_t(F))$ is free abelian with the basis
$\bar{Y_t}=\{wL_t\gamma_2(\rho_t(F))\mid \ w\in E_t\}.$ By
considering the abelian group ${\cal N}_{c_1,c_2,\ldots,c_t}M(G)$
and the form of the elements of $L_t$, we have
$L_t\gamma_2(\rho_t(F))\leq [R,\ _{c_1}F,\
_{c_2}\rho_1(F),\ldots,\ _{c_t}\rho_{t-1}(F)],$ and the following
isomorphism
$${\cal N}_{c_1,c_2,\ldots,c_t}M(G)\cong \frac{\rho_t(F)/L_t\gamma_2(\rho_t(F))}{[R,\ _{c_1}F,\ _{c_2}\rho_1(F),\ldots,\ _{c_t}\rho_{t-1}(F)]/L_t\gamma_2(\rho_t(F))}.$$

Clearly $\gamma_{\pi}(F)/\gamma_{\pi+1}(F)$ is the free abelian
group with the basis of all basic commutators of weight $\pi$ on
$X$, where $\pi=\prod_{i=1}^{t} (c_i+1)$. Using the above fact and
$\bar{Y_t}$, the basis of the free abelian group
$\rho_t(F)/L_t\gamma_2(\rho_t(F))$, one can obtain the following
inclusion:
$$\gamma_{\pi+1}(F)\cap \rho_t(F)\leq L_t\gamma_2(\rho_t(F)).$$
Therefore, it is clear that $\bar{Y_t}=C_{0,t}\cup C_{1,t}\cup
\ldots \cup C_{k,t}$ is a basis for the free abelian group
$\rho_t(F)/L_t\gamma_2(\rho_t(F))$ and $D_{0,t}\cup D_{1,t}\cup
\ldots \cup D_{k,t}$ is a basis for the free abelian group
$[R,_{c_1}F,_{c_2}\rho_1(F),\ldots,_{c_t}\rho_{t-1}(F)]/L_t\gamma_2(\rho_t(F)),$
where
\begin{center}
$C_{0,t}=\{wL_t\gamma_2(\rho_t(F))\in \bar{Y_t}\mid w\in E_t $ and
$w$ is a\\
commutator on letters $x_1,\ldots,x_m\};$\\
 $C_{\lambda,t}=\{wL_t\gamma_2(\rho_t(F))\in \bar{Y_t}\mid \ w\in
E_t$ and $w$ is a commutator on letters
$x_1,\ldots,x_m,x_{m+1},\ldots,x_{m+\lambda}$
 such that the letter $x_{m+\lambda}$ does appear in $w\};$

$D_{\lambda,t}=\{w^{n_{\lambda}}L_t\gamma_2(\rho_t(F))\mid \
wL_t\gamma_2(\rho_2(F))\in C_{\lambda,t}\};$\\ for all  $1\leq
\lambda \leq k.$
\end{center}
Note that $\mid C_{0,t}\mid
=\chi_{c_t+1}(\ldots(\chi_{c_1+1}(m))\ldots)$ and

$\mid C_{\lambda,t}\mid
=\chi_{c_t+1}(\ldots(\chi_{c_1+1}(m+\lambda))\ldots)-\chi_{c_t+1}(\ldots(\chi_{c_1+1}(m+\lambda-1))\ldots).$
Hence the result holds.$\Box$

Now we can state the following interesting corollary.\\

\hspace{-0.6cm}{\bf Corollary 2.4.} Let ${\cal S}_{\ell}$ be the
variety of solvable groups of length at most $\ell$ and $G\cong
{\bf Z}^{(m)}\oplus {\bf Z}_{n_1}\oplus {\bf Z}_{n_2}\oplus \ldots
,{\bf Z}_{n_k}$ be a finitely generated abelian group, where
$n_{i+1}\mid n_i$ for all $1\leq i \leq k-1$. Then the following
isomorphism holds:

$${\cal S}_{\ell}M(G)\cong {\bf Z}^{(h_m)}\oplus {\bf Z}_{n_1}^{(h_{m+1}-h_m)}\oplus\ldots\oplus
{\bf Z}_{n_k}^{(h_{m+k}-h_{m+k-1})}$$

where $h_i=\chi_2\underbrace{(\ldots(
}_{(l-times)}\chi_2(i))\ldots)$ for all $m\leq i \leq m+k.$

\end{document}